\documentclass[a4paper,fleqn,10pt,twocolumn]{SICE_ISCS}
\usepackage{amsmath}
\usepackage{amsfonts}
\usepackage[numbers]{natbib}
\usepackage{xcolor}

\title{Equivalence of additive and parametric pinning control protocols for systems of weakly coupled oscillators}

\author{Riccardo Muolo${}^{1,2\dagger}$ and Yuzuru Kato${}^{3}$}
\speaker{Riccardo Muolo}

\affils{${}^{1}$RIKEN Center for Interdisciplinary Theoretical and Mathematical Sciences (iTHEMS), Saitama, Japan\\
${}^{2}$Department of Systems and Control Engineering, Institute of Science Tokyo, Tokyo, Japan\\
(E-mail: riccardo.muolo@riken.jp)\\
${}^{3}$Department of Complex and Intelligent Systems, Future University Hakodate, Hokkaido, Japan\\
(E-mail: katoyuzu@fun.ac.jp)\\
}
\abstract{%
Controlling the behavior of nonlinear systems on networks is a paramount task in control theory, in particular the control of synchronization, given its vast applicability. In this work, we focus on pinning control and we examine two different approaches: the first, more common in engineering applications, where the control is implemented through an external input (additive pinning); the other, where the parameters of the pinned nodes are varied (parametric pinning). By means of the phase reduction technique, we show that the two pinning approaches are equivalent for weakly coupled systems exhibiting periodic oscillatory behaviors. Through numerical simulations, we validate the claim for a system of coupled Stuart--Landau (SL) oscillators. Our results pave the way for further applications of pinning control in real-world systems.}


\begin{document}

\maketitle

%-----------------------------------------------------------------------

\section{Introduction}

Synchronization phenomena in networks of coupled oscillators have been widely studied in complex systems, due to their importance in fields such as physics, biology, and engineering~\cite{pikovsky2001synchronization,arenas2008synchronization}. Particularly relevant is the control of synchronization~\cite{liu2016control}, which is a challenging task due to the nonlinear nature of oscillatory systems, and complex systems in general~\cite{d2023controlling,coraggio2026controlling}.
A powerful approach in this context is \emph{pinning control}, which consists in applying an external input to a selected subset of nodes in the network in order to drive the dynamics of the entire ensemble towards a desired (i.e., target) state \cite{grigoriev1997pinning,sorrentino2007controllability}. Pinning control has been successfully applied in the framework of opinion dynamics \cite{iudice2022bounded,ancona2023influencing}, epidemics \cite{du2015selective,yang2019feedback}, pattern formation \cite{buscarino2019turing}, synchronization dynamics \cite{porfiri2008criteria,yu2013synchronization}, and chimera states \cite{gambuzza2016pinning}. Recently, such technique has been applied also in the framework of higher-order interactions~\cite{battiston2021physics,millan2025topology}, in particular to control synchronization dynamics~\cite{della2023emergence,de2022pinning,de2023pinning,muolo2025pinning}. In the latter work~\cite{muolo2025pinning}, two distinct protocols have been considered: \emph{additive pinning}, where a control signal is directly added to the dynamical variables of the pinned nodes, and \emph{parametric pinning}, where the control acts by modifying intrinsic parameters, such as the natural frequency, of the pinned nodes. Note that additive pinning is the usual pinning control found in the literature, while the parametric pinning is less used due to the difficulties in the implementations. However, in social systems, e.g., opinion dynamics or epidemics, it may be easier to vary the parameters (e.g., introducing some control from within, which acts on the parameters of the system) rather than apply an external forcing. 

In this work, we study networks of Stuart--Landau (SL) oscillators, which provide the normal form for limit-cycle dynamics near a Hopf-Andronov bifurcation \cite{Kuramoto,nakao2014complex}. For weakly coupled systems, the {multi-dimensional} dynamics of each oscillator can be reduced to a {single-dimensional}
phase equation using the phase reduction method~\cite{nakao2016phase,monga2019phase}. In this approximation, each oscillator $i$ is represented solely by its phase $\vartheta_i \in [0,2\pi)$. The phase reduction allows us to analyze the dynamics in a lower-dimensional framework while preserving essential features of the synchronization behavior, as long as the coupling and the perturbations involved are weak. 

Using this framework, we implement both additive and parametric pinning control and show that the two protocols are equivalent for the phase model. Moreover, we numerically show that, when the validity conditions of the phase reduction are satisfied (i.e., weak coupling and weak perturbations), the two protocols produce equivalent results also for the non-reduced models. The method developed has not only a theoretical interest, but can prove useful in applications when the control of synchronization is important.

In the next Section we introduce the model, while in Sec. \ref{sec:phase_red} we discuss the theory of phase reduction. After having introduced the pinning protocols in Sec. \ref{sec:pinning}, we show that they are equivalent for phase models in Sec. \ref{sec:equiv}, where we also numerically show the equivalence also for non-reduced model, right before concluding and laying down possible future directions of investigation.

\section{Networks of coupled Stuart--Landau (SL) oscillators}\label{sec:Stuart--Landau}

We consider a system of $n$ identical interacting SL units. The dynamics of each uncoupled oscillator are given by
\begin{align}\label{eq:Stuart--Landau_isolated}
    \dot{x}_i &= f(x_i,y_i) = \alpha x_i - \omega y_i - (x_i^2 + y_i^2)x_i, \nonumber \\
    \dot{y}_i &= g(x_i,y_i) = \omega x_i + \alpha y_i - (x_i^2 + y_i^2)y_i,
\end{align}
where $\alpha$ is the bifurcation parameter (for $\alpha>0$ a  stable limit cycle of amplitude $\sqrt{\alpha}$ exists) and $\omega$ is the natural frequency {of the oscillator}. Note that all results developed in this work remain valid for small heterogeneities across the nodes, i.e., $\omega_i=\omega+\delta\omega_i$ and $\alpha_i=\alpha+\delta\alpha_i$, with $\delta\omega_i,\delta\alpha_i<<1$ for all $i$. By denoting $\vec{X}_i=(x_i,y_i)^\top$ and $\vec{F}=(f,g)^\top$, we can rewrite the above equation in vector form \begin{displaymath}
    \dot{\vec{X}}_i=\vec{F}(\vec{X}_i).
\end{displaymath}

We describe the network structure by a (weighted) adjacency matrix $A=\{A_{ij}\}$, where $A_{ij}\ge0$ denotes the strength of the symmetric interaction between nodes $j$ and $i$. The degree of node $i$ is defined as
\begin{displaymath}
  k_i = \sum_{j=1}^n A_{ij}.
\end{displaymath}
We can now introduce the (negative semidefinite) Laplacian matrix $L$ by
\[
  L_{ij}=\begin{cases}
    A_{ij}, & i\neq j, \\
    -k_i, & j=i,
  \end{cases}
\]
so that $L\,\vec{1}=0$, where $\vec{1}$ is an $n$-dimensional homogeneous vector with all entries equal to~$1$.

We consider diffusive diagonal coupling. Denoting by $\varepsilon>0$ the coupling strength, the equations for a network of coupled SL oscillators are given by
\begin{equation}\label{eq:Stuart--Landau_network_new}
\begin{aligned}
\dot{x}_i &= f(x_i,y_i) + \varepsilon \sum_{j=1}^n L_{ij}x_j, \\[4pt]
\dot{y}_i &= g(x_i,y_i) + \varepsilon \sum_{j=1}^n L_{ij}y_j.
\end{aligned}
\end{equation}

In vector form, system \eqref{eq:Stuart--Landau_network_new} becomes
\begin{equation}\label{eq:vector_form_new}
\dot{\vec{X}}_i = \vec{F}(\vec{X}_i) + {\mathbf{D}}\sum_{j=1}^n L_{ij}\vec{X}_j,
\end{equation}

\noindent where $\mathbf{D}$ is the diagonal coupling matrix given by \begin{displaymath}
    \mathbf{D}=\varepsilon\begin{bmatrix}
        1 & 0 \\ 0 & 1
    \end{bmatrix}.
\end{displaymath}
Note that the coupling matrix $\mathbf{D}$ can take different forms~\cite{nakao2014complex,muolo2024phase}. In the context of the equivalence between additive and parametric pinning protocols, the form of $\mathbf{D}$ does not affect the results. Hence, we will consider the form of the above equation. 

By introducing diffusive coupling functions\\ $\vec{H}_{ij}=\vec{H}(\vec{X}_j,\vec{X}_i)=\vec{X}_j-\vec{X}_i$, we can rewrite Eq. \eqref{eq:vector_form_new} in a more general form, which will turn useful in the next Section:

\begin{equation}
    \dot{\vec{X}}_i = \vec{F}(\vec{X}_i) + {\mathbf{D}}\sum_{j=1}^n A_{ij}\vec{H}(\vec{X}_j(\vartheta_j), \vec{X}_i(\vartheta_i)),
\end{equation}

\section{The Phase Reduction Approach}\label{sec:phase_red}

We will now briefly introduce the theory of phase reduction. For more details, see {Refs.} \cite{Kuramoto,nakao2014complex,monga2019phase}. 

Let us consider a network of $n$ coupled oscillators as in the previous section.
Under the assumption of weak coupling ($\varepsilon \ll 1$) and that each oscillator exhibits a stable limit cycle, we can perform the phase reduction. Let $\vartheta_i \in [0,2\pi)$ denote the phase of the $i$th oscillator along its limit cycle. Our goal is to obtain the following phase-reduced dynamics
\begin{equation}\label{eq:phase_network}
    \dot{\vartheta}_i = \omega + \varepsilon \sum_{j=1}^n A_{ij} \vec{Z}(\vartheta_i) \cdot \vec{\mathcal{H}}(\vec{X}_j(\vartheta_j), \vec{X}_i(\vartheta_i)),
\end{equation}
where $\vec{\mathcal{H}}$ is a phase coupling function between $i$th and $j$th oscillator,
$\vec{Z}(\vartheta_i) \in \mathbf{R}^2$ is the phase sensitivity function (PSF) of the $i$th oscillator whose phase $\vartheta_i$, which determines whether a perturbation slows down or accelerates the phase of the oscillations in each point of the limit cycle, and $\vec{X}_i(\vartheta_i)$ represents the state of the $i$th oscillator on the limit cycle. Note that, since in our case the oscillators are identical, the PSF $\vec{Z}(\vartheta_i)$ is the same for every oscillator $i$. For SL oscillators of Eq. \eqref{eq:Stuart--Landau_isolated}, the PSF can be explicitly computed~\cite{nakao2016phase} and takes the form
\begin{equation}
    \vec{Z}(\vartheta) = \frac{1}{\sqrt{\alpha}}(-\sin\vartheta, \cos\vartheta)^\top.
\end{equation}
In general, the PSF can be obtained only numerically, for example using the adjoint method~\cite{nakao2016phase}.

More formally, given the following system 
\begin{displaymath}
    \dot{\vec{X}}_i = \vec{F}(\vec{X}_i) + {\mathbf{D}}\sum_{j=1}^n A_{ij}\vec{H}(\vec{X}_j,\vec{X}_i),
\end{displaymath}
assume that: each uncoupled oscillator ($\varepsilon = 0$) has a stable limit cycle with period $T$ and phase $\vartheta_i \in [0,2\pi)$, the coupling is sufficiently weak, so that the trajectories remain close to the individual limit cycles.

Then, to first order in $\varepsilon$, the high-dimensional system can be approximated by the following phase-reduced equations
\begin{equation}
    \dot{\vartheta}_i = \omega + \sum_{j=1}^n A_{ij} \vec{Z}(\vartheta_i) \cdot \vec{\mathcal{H}}(\vec{X}_j(\vartheta_j), \vec{X}_i(\vartheta_i)),
\end{equation}
where $\vec{\mathcal{H}}=\boldsymbol{D}\cdot\vec{H}(\vec{X}_j(\vartheta_j), \vec{X}_i(\vartheta_i))$, $\omega$ is the natural frequency of the oscillators, and $\vec{Z}(\vartheta_i)$ is the phase sensitivity function of the $i$th oscillator, satisfying $\vec{Z}(\vartheta_i) \cdot \vec{F}(\vec{X}(\vartheta_i)) = \omega$.

It is worth noting that, in our case of SL oscillators with linear diffusive coupling, i.e., $\vec{H}(\vec{X}_j,\vec{X}_i)=\vec{X}_j-\vec{X}_i$, the phase model obtained is exactly the Kuramoto model~\cite{kuramoto1975}, namely

\begin{equation}\label{eq:kura_saka}
    \dot{\vartheta}_i = \omega + \varepsilon \sum_{j=1}^n A_{ij} \sin(\vartheta_j-\vartheta_i).
\end{equation}
Note that there is no phase-lag because of the chosen form of the SL. In fact, with other coupling configurations, one obtains the Kuramoto-Sakaguchi model~\cite{sakaguchi1986soluble}, i.e., the Kuramoto model with a phase lag.

Let us point out that, in general, after having obtained the PSF and computed the scalar product $\vec{Z}(\vartheta_i) \cdot \vec{\mathcal{H}}(\vec{X}_j(\vartheta_j), \vec{X}_i(\vartheta_i))$, {we often}
perform an averaging procedure~\cite{averaging_sanders} in order to remove non-resonant terms.

Lastly, let us note that the theory of phase reduction has been recently extended also to the framework of higher-order interactions~\cite{leon2024higher,leon2025phase}, meaning that all results that we obtain for oscillators coupled via a network remain valid when the coupling is higher-order, i.e., via a hypergraph or a simplicial complex.

\section{Additive and Parametric Pinning Control Protocols}\label{sec:pinning}

We use two control schemes to implement pinning in the network: \emph{additive pinning} and \emph{parametric pinning}, following~Ref.~\cite{muolo2025pinning}. Both rely on two shared parameters: the set of pinned nodes and the time window of control. We denote by $t_p$ the duration of the control and by $n_p<n$ the number of pinned nodes. The set of all nodes is $I=\{1,\ldots,n\}$, and the pinned nodes are $I_p$. Unless otherwise stated, we take $I_p=\{1,\ldots,n_p\}$.

\subsection{Additive pinning}

This method, used for example in~\cite{gambuzza2016pinning}, applies an external input to selected nodes. The network equations become
\begin{equation}\label{eq:additive_network}
    \dot{\vec{X}}_i=\vec{F}(\vec{X}_i)+{\mathbf{D}}\sum_{j=1}^n L_{ij}\,\vec{X}_j +\vec{U}_i,
\end{equation} 
where $L$ is the Laplacian matrix and $\vec{U}_i$ is the control input:
\begin{equation}
{\vec{U}_i=\left \{ 
    \begin{array}{l}
    \vec{0}, \quad i \in I\setminus I_p,  \\
    \vec{u}_i, \quad i \in I_p,
    \end{array}
    \right.}
\end{equation}
with $\vec{0}=(0,0)^\top$, $\vec{u}_i=(u_{i}(t),u_{i}(t))^\top$, and
\begin{equation}
    u_{i}(t)= \lambda_{i}[\Theta(t)-\Theta(t-t_p)],
\end{equation}
where $\lambda_{i}$ are drawn from a uniform distribution and $\Theta$ is the Heaviside step function.

For pinned nodes $i\in I_p$, the controlled SL dynamics reads
\begin{equation}
\begin{cases}
\displaystyle \dot{x}_i =f + \varepsilon \sum_{j=1}^n L_{ij}x_j+\lambda_{i}[\Theta(t)-\Theta(t-t_p)],\\
\displaystyle \dot{y}_i =g + \varepsilon \sum_{j=1}^n L_{ij}y_j+\lambda_{i}[\Theta(t)-\Theta(t-t_p)],
\end{cases}
\label{eq:network_controlledI}
\end{equation}
while unpinned nodes ($i\in I\setminus I_p$) follow the original dynamics.

\subsection{Parametric pinning}

In this scheme, the control acts by modifying parameters of the pinned nodes. The network equations are
\begin{equation}\label{eq:parametric_network}
\dot{\vec{X}}_i=\vec{F}_i(\vec{X}_i)+{\mathbf{D}}\sum_{j=1}^n L_{ij}\vec{X}_j,
\end{equation}
where
\begin{equation}
{\vec{F}_i=\left \{ 
    \begin{array}{l}
    \vec{F}, \quad i \in I\setminus I_p,  \\
    (f_{p},g_{p})^\top, \quad i \in I_p.
    \end{array}
    \right.}
\end{equation}

The dynamics of pinned nodes follow
\begin{equation}\label{eq:Stuart--Landau_pinned_network}
\begin{cases}
 f_{p}( x_{i}, y_{i}) = \alpha  x_{i} - \Omega_{p,i} y_{i}  -\left(x_{i}^2  + y_{i}^2 \right)x_{i}, \\\\
 g_{p}( x_{i}, y_{i}) = \Omega_{p,i} x_{i} + \alpha y_{i} - \left(x_{i}^2  + y_{i}^2 \right)y_{i},
\end{cases}
\end{equation}
whose frequencies are given by
\begin{equation}\label{eq:new_omega_network}
    \Omega_{p,i}(t)=\begin{cases}
    \omega_{p,i}, & t\leq t_p,  \\
    \omega, & t>t_p,
    \end{cases}
\end{equation} 
where $\omega_{p,i}$ are drawn from a uniform distribution. During the time in which the control is active, i.e., $t\leq t_p$, the frequency is $\omega_{p,i}$, and returns to $\omega$ afterwards. For unpinned nodes, the frequency is always $\omega$.

\section{Equivalence of the two pinning protocols}\label{sec:equiv}

In this Section, we will show the equivalence of the two pinning protocols when the framework is such that the phase reduction approximation is valid, i.e., weak coupling strengths and weak perturbations. We will first show that the two protocols are exactly the same when considering phase models, to then prove that they are equivalent also for non-phase models, as long as the phase reduction approximation is valid.

\subsection{Formal equivalence for phase models}

Let us consider the Kuramoto model of coupled phase oscillators~\cite{kuramoto1975}
\begin{displaymath}
    \dot{\vartheta}_i = \omega + \varepsilon \sum_{j=1}^n A_{ij} \sin(\vartheta_j-\vartheta_i).
\end{displaymath}

If we apply additive pinning according to the first protocol, we obtain the following dynamics for pinned nodes

\begin{displaymath}
    \dot{\vartheta}_i = \omega + \varepsilon \sum_{j=1}^n A_{ij} \sin(\vartheta_j-\vartheta_i)+\lambda_{i}[\Theta(t)-\Theta(t-t_p)].
\end{displaymath} 

If then apply parametric pinning according to the second protocol, we obtain that, for $t\leq t_p$, pinned nodes follow 

\begin{displaymath}
    \dot{\vartheta}_i = \omega_{p,i} + \varepsilon \sum_{j=1}^n A_{ij} \sin(\vartheta_j-\vartheta_i).
\end{displaymath}

By setting $\omega_{p,i}=\omega+\lambda_i$, it is clear that the two protocols are exactly the same. In fact, the dynamics being linear (i.e., given by the sole frequency), the control input is linearly additive to the dynamics. Note that the same reasoning can be applied to the case of heterogeneous frequencies, i.e., $\omega_i$ for the $i$th oscillator.

Moreover, this remains valid for any phase model. In fact, even when the PSF cannot be computed analytically as for the SL, the phase model will still be such that each isolated oscillator will follow its intrinsic frequency, which is a key assumption of the phase reduction. The overall model, called, in this case, Winfree model~\cite{Win80} or Kuramoto-like model, will have a more complicated interaction term, in general not treatable analytically.

Lastly, let us observe that, despite the equivalency between phase models, independently of the interaction terms, we should be aware that the phase models are an approximation of real-world oscillatory systems, which do not completely follow a Kuramoto-like dynamics, but are rather approximated by it. In fact, the phase reduction serves the purpose of characterizing the dynamics of oscillatory systems, offering a very good approximation when the coupling and the perturbations involved are weak. However, the equivalence of the two pinning approaches must be proven for non-reduced systems to be considered applicable to oscillatory systems. This is the task we will carry out in what follows.

\subsection{Numerical equivalence for non-reduced models}

Let us consider a system of $n$ Stuart-Landau oscillators coupled via a network, and apply additive pinning control. The dynamics is described by Eqs. \eqref{eq:network_controlledI}. If we perform the phase reduction, we obtain the following phase equation for $t<t_p$

\begin{equation}\label{eq:pase_pinn}
     \dot{\vartheta}_i = \omega + \varepsilon \sum_{j=1}^n A_{ij} \sin(\vartheta_j-\vartheta_i)+\sqrt{\frac{2}{\alpha}}\lambda_i\cos\Big(\vartheta + \frac{\pi}{4}\Big).
\end{equation}

From the above equation, we see that the frequency $\omega$ of each phase-reduced oscillator $i$ receives an input of $\sqrt{\frac{2}{\alpha}}\lambda_i\cos\Big(\vartheta + \frac{\pi}{4}\Big)$ until $t=t_p$ and the pinning control is stopped. Hence, if we set the frequency of the parametric pinning given by Eq. \eqref{eq:new_omega_network} to \begin{equation}
    \omega_{p,i}=\omega+\frac{1}{t_p}\int_{0}^{t_p} \sqrt{\frac{2}{\alpha}}\lambda_i\cos\Big(\vartheta(t) + \frac{\pi}{4}\Big) dt,
\end{equation}
we have that the two pinning approaches are equivalent in the limit of validity of the phase reduction, i.e., $\varepsilon<<1$ and $\lambda_i<<1$.

\begin{figure*}[t!]
    \centering
    \includegraphics[width=\textwidth]{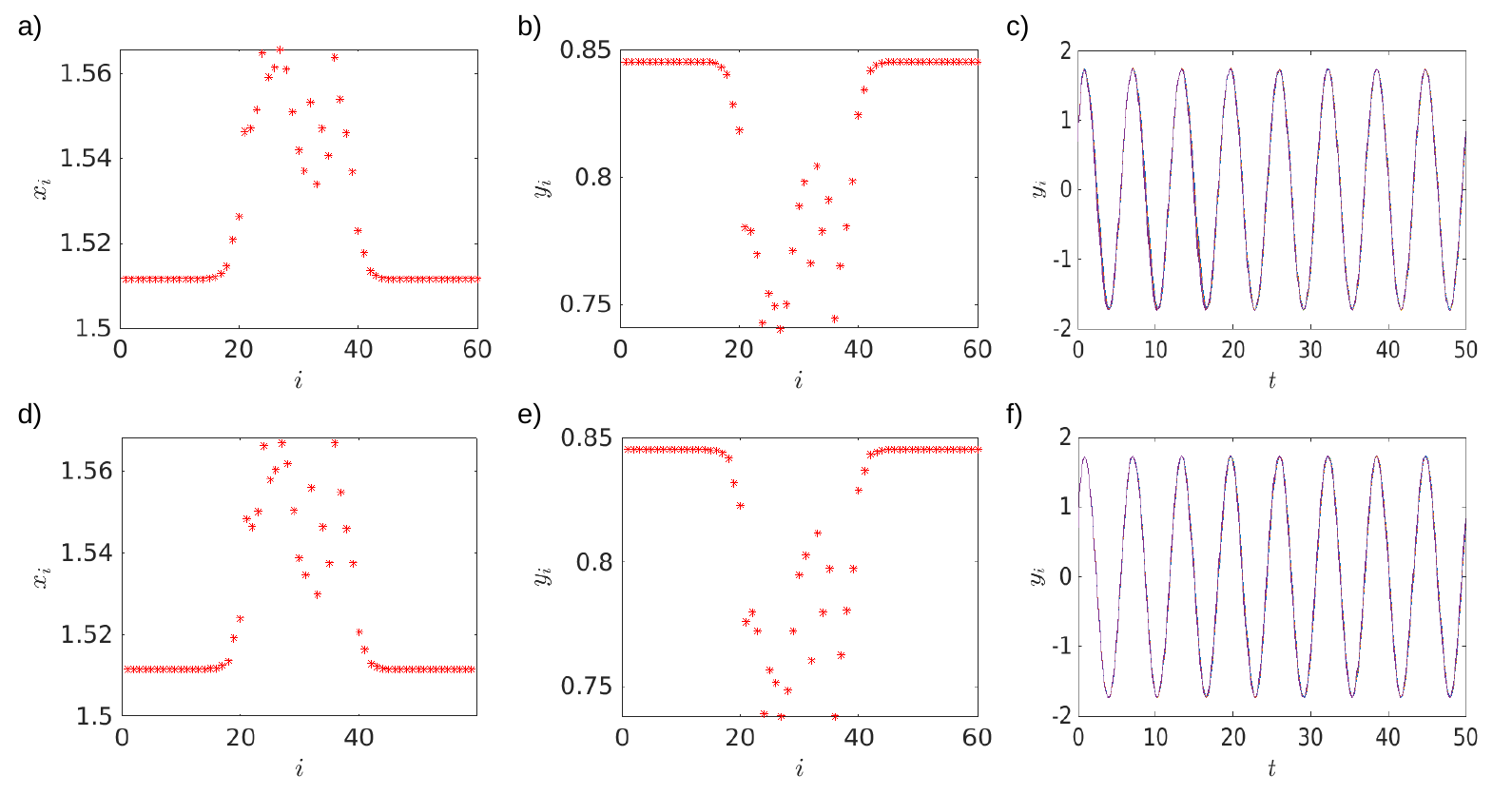}
    \caption{In this setting, the phase reduction approach is valid and so is the equivalence of the two pinning protocols. The upper panels show the case of additive pinning, while the lower panels show the parametric pinning. Panels a) and d) show a snapshot of the $x_i$ variables for $50$ time units (t.u.); panels b) and e) show a snapshot of the $y_i$ variables for $50$ time units (t.u.); panels c) and f) show the time series of the $y_i$ variables. The network is a $1$-dimensional $4$-regular lattice of $60$ nodes; the parameters of the SL {oscillators} are $\alpha=3$ and $\omega=1$; the coupling is $\varepsilon=0.01$ and the order of the pinning perturbation is $\lambda=0.1$. The system is integrated for $50$ t.u. with the Runge-Kutta IV method with integration step $dt=0.01$; the pinning is applied to $n_p=20$ nodes for $t_p=10$ t.u.}
    \label{fig:equiv}
\end{figure*}

\begin{figure*}[t!]
    \centering
    \includegraphics[width=\textwidth]{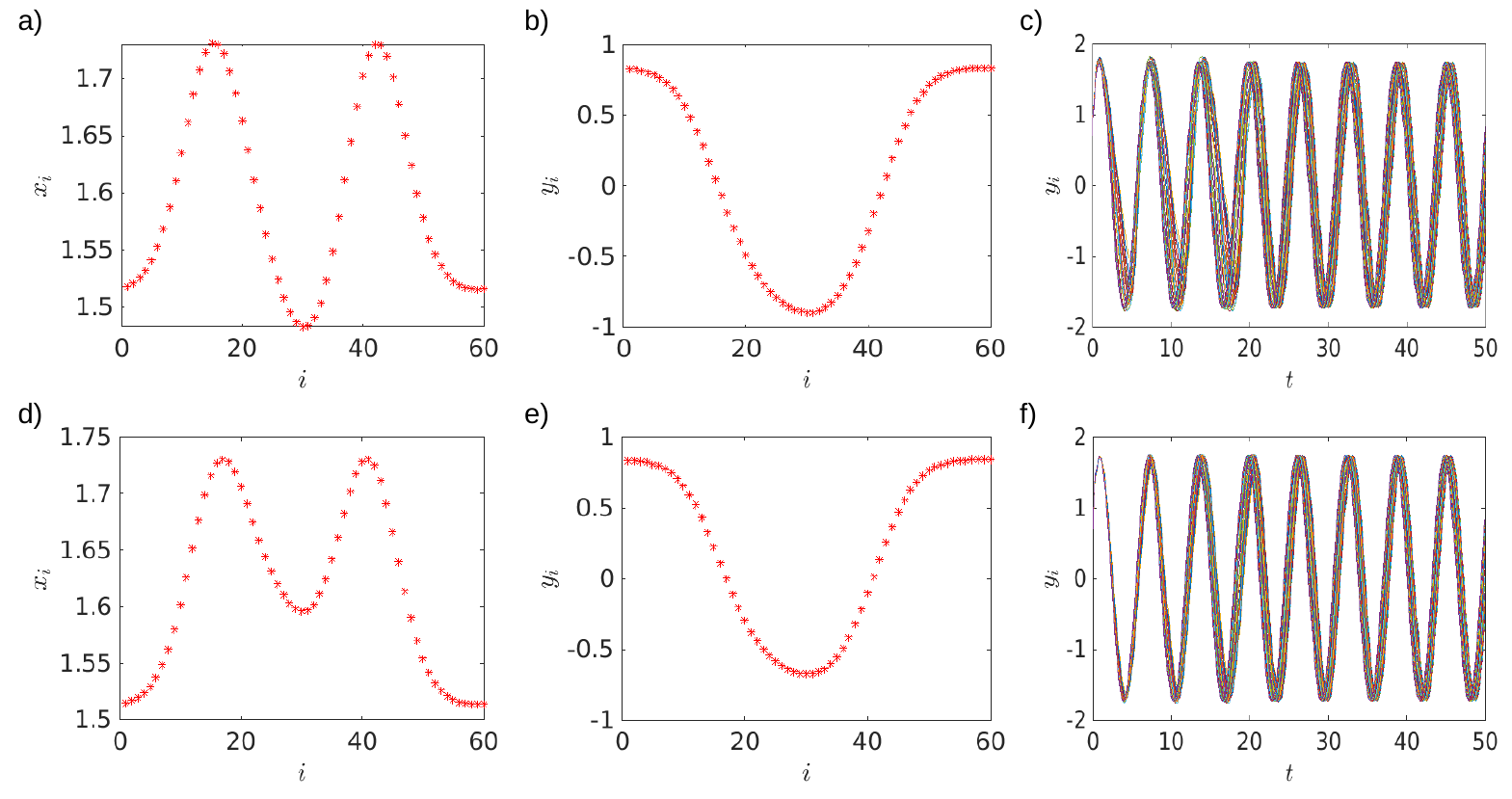}
    \caption{In this setting, the phase reduction approach is not anymore valid and, hence, the two pinning protocols are not equivalent. The upper panels show the case of additive pinning, while the lower panels show the parametric pinning. Panels a) and d) show a snapshot of the $x_i$ variables for $50$ time units (t.u.); panels b) and e) show a snapshot of the $y_i$ variables for $50$ time units (t.u.); panels c) and f) show the time series of the $y_i$ variables. The network is a $1$-dimensional $4$-regular lattice of $60$ nodes; the parameters of the SL {oscillators} are $\alpha=3$ and $\omega=1$; the coupling is $\varepsilon=0.2$ and the order of the pinning perturbation is $\lambda=0.6$. The system is integrated for $50$ t.u. with the Runge-Kutta IV method with integration step $dt=0.01$; the pinning is applied to $n_p=20$ nodes for $t_p=10$ t.u.}
    \label{fig:nonequiv}
\end{figure*}

We can now validate the above result through simulations of SL oscillators coupled via a network, namely, a regular lattice. In our numerical study, we examine two different settings: one in which the assumptions of validity of the phase reduction are satisfied, whose results are shown in Fig. \ref{fig:equiv}, and the other, whose results are shown in Fig. \ref{fig:nonequiv}, in which the phase reduction approximation is not anymore valid.
Fig. \ref{fig:equiv} shows that the two approaches are almost indistinguishable when phase reduction is valid. On the other hand, when the phase reduction is not valid, the two pinning protocols give different results, as shown in Fig. \ref{fig:nonequiv}. 

This phenomenology, which was expected, can be intuitively understood when looking at the time series of panels c) and f) of the two Figures. In the case where the phase reduction is valid, the amplitude remains almost constant, hence, the phase model fully describes the dynamics. However, when the coupling and the perturbations are not weak, the additive pinning affects the amplitude, while the parametric pinning does not, hence, the phase models cannot give an accurate description of the dynamics.  

Let us add that this result is valid in general for any oscillatory system undergoing a supercritical Hopf-Andronov bifurcation, because the SL oscillator is the normal form of such bifurcation~\cite{nakao2014complex}. This means that every such oscillator, e.g., the FitzHugh-Nagumo model, the van der Pol equation, etc., behave as a SL close to the Hopf-Andronov bifurcation. Additionally, phase reduction is valid for any limit cycle oscillator as long as the coupling and the perturbations are weak, making this framework valid for general oscillators, not only close to the Hopf-Andronov bifurcation. Lastly, note that we have shown the case of regular lattice for {the} sake of simplicity, but the results remain valid for non-regular networks and also for hypergraphs and simplicial complexes, because the network topology has no role in the method. For all these reasons, we can conclude that our results are more general than what we have presented in this section and still stand for different oscillatory systems coupled via different network or higher-order topologies.

\section{Conclusion}

In this work, we have shown the equivalence of two different pinning control protocols, additive pinning and parametric pinning, for oscillatory and periodic weakly coupled systems. While additive pinning is often implemented in engineering applications, parametric pinning seems more difficult {to implement} in practice. In fact, in experiments, it is easier to add an external input, while the task of acting on the parameters of the systems to be controlled seems more challenging. Nonetheless, in certain systems, such as social systems, it is easier to directly act on the model parameters rather than applying an external input. For instance, in case of an epidemic spreading, one can control a given location, i.e., a node of the network, by implementing certain containment measures. In practice, this would mean acting on the parameters of that given node, i.e., parametric pinning. It is true that we have shown the equivalence of the two protocols only for oscillatory systems and not in general. However, oscillatory behaviors occur also in such social systems and epidemics: for instance, when considering a periodic infectivity parameter, it was shown that the system exhibits a periodic attractor~\cite{boatto2018sir}. Therefore, our result is twofold: from a purely theoretical point of view, it deepens our knowledge of the phase reduction method for systems to which a control is applied; on the other hand, it also provides guidance for applications, by formally proving the equivalence of the two control protocols, once certain conditions are satisfied.

Future studies could further explore such framework and implement parametric pinning control to contain the epidemics. 
Lastly, the framework hereby studied could be further extended to consider non-oscillatory or non-periodic systems, and investigate under which conditions the two protocols remain equivalent.

\section*{Acknowledgements}
The authors are grateful to Hiroya Nakao, Mattia Frasca, and Lucia Valentina Gambuzza for discussions. R.M. acknowledges JSPS KAKENHI 24KF0211 for financial support. Y.K. acknowledges JSPS KAKENHI JP22K14274, and JST PRESTO JPMJPR24K3 for financial support.

%%%%%%%%%%%%%%%%% BIBLIOGRAPHY IN THE LaTeX file !!!!! %%%%%%%%%%%%%%%%%%%%%%
%\bibliography{biblio}

\end{document}